\documentclass[a4paper,11pt]{article}

\usepackage{amsmath,amsthm,amssymb,mathrsfs,latexsym,amsfonts}
\usepackage{graphicx,psfrag,epsfig}
\usepackage[english]{babel}
\usepackage[latin1]{inputenc}

\newtheorem{theoreme}{Theorem}[section]

\newtheorem{proposition}[theoreme]{Proposition}

\newtheorem{definition}[theoreme]{Definition\rm}
\newtheorem{remarque}[theoreme]{\it Remark}

\begin{document}

\bibliographystyle{amsalpha}

\def\MP{\,{<\hspace{-.5em}\cdot}\,}
\def\SP{\,{>\hspace{-.3em}\cdot}\,}
\def\PM{\,{\cdot\hspace{-.3em}<}\,}
\def\PS{\,{\cdot\hspace{-.3em}>}\,}
\def\EP{\,{=\hspace{-.2em}\cdot}\,}
\def\PP{\,{+\hspace{-.1em}\cdot}\,}
\def\PE{\,{\cdot\hspace{-.2em}=}\,}
\def\N{\mathbb N}
\def\C{\mathbb C}
\def\Q{\mathbb Q}
\def\R{\mathbb R}
\def\T{\mathbb T}
\def\A{\mathbb A}
\def\Z{\mathbb Z}
\def\demi{\frac{1}{2}}
%\phantom{}
%\vskip-10truecm
\begin{titlepage}
\author{Laurent Niederman~\footnote{Laurent.Niederman@math.u-psud.fr, Topologie et Dynamique, Laboratoire de Math\'ematiques d'Orsay (UMR-CNRS 8628), Universit\'e Paris Sud, 91405 Orsay Cedex.}
 {}~\footnote{Astronomie et Syst\`emes Dynamiques, Institut de M\'ecanique C\'eleste et de calcul des \'eph\'em\'erides (UMR-CNRS 8028), Observatoire de Paris, 75014
Paris.}}
\title{\LARGE{\textbf{Effective generic super-exponential stability of elliptic
fixed points for symplectic vector fields.}}}
\vskip3truecm
\end{titlepage}

\maketitle

\begin{abstract}
In this article, we consider linearly stable elliptic fixed points (equilibrium) for a symplectic vector field and we prove generic results of super-exponential stability for nearby solutions. We will focus on the neighbourhood of elliptic fixed points but the case of linearly stable isotropic reducible invariant
tori in an Hamiltonian system should be similar.

More specifically, Morbidelli and Giorgilli have proved a result of stability over super-exponentially long times if one consider an analytic lagrangian torus, invariant for an analytic hamiltonian system, with a diophantine translation vector which admit a sign definite torsion. Then, the solutions of the system move very little over times which are super-exponentially long with respect to the inverse of the distance to the invariant torus.

 The proof is in two steps: first the construction of a Birkhoff normal form at a high order, then the application of Nekhoroshev theory. Bounemoura has shown that the second step of this construction remains valid if the Birkhoff normal form linked to the invariant torus or the elliptic fixed point belongs
to a generic set among the formal series. This is not sufficient to prove this kind of super-exponential stability results in a general setting. We should also establish that most strongly non resonant elliptic fixed point or invariant torus in a Hamiltonian system admit a Birkhoff normal form in the set introduced by Bounemoura. We show here that this property is satisfied generically in the sense of the measure (prevalence) through infinite-dimensional probe spaces (that is an infinite number of parameter chosen at random) with methods similar to those developed in a paper of Gorodetski, Kaloshin and
Hunt in another setting.

\end{abstract}

\section{Introduction}

  We are interested in the stability properties, in the sense of Lyapounov, of linearly stable elliptic fixed points (equilibrium) denoted
$x_*$ on a symplectic manifold $(M,\Omega)$ for a symplectic vector field $X$ (i.e. : $i_X\Omega$ is closed), that is $X(x_*)=0$.

As the problem is local, it is enough to consider a Hamiltonian $H$, defined and analytic on an open neighbourhood of $0 \in\R^n\times\R^n$ equipped with
the standard symplectic form, having the origin as a fixed point. Up to an irrelevant additive constant, expanding the Hamiltonian as a power series at the origin we can write
\[ H(z)=H_2(z)+ V(z) \]
where $z=(x,y)\in\R^n\times\R^n$ is sufficiently close to the origin in $\R^n\times\R^n$, $H_2$ is the Hessian of $H$ at $0$ and $V(z)=O(|z|^3)$.

 Recall that the fixed point is said to be elliptic if the spectrum of the linearized system is purely imaginary. Due to the symplectic character of the equations, these are the only linearly stable fixed points and the spectrum has the form $\{\pm i\omega_1, \dots, \pm i\omega_n\}$ for some vector $\omega=(\omega_1, \dots, \omega_n)$ which is called the normal (or characteristic) frequency. Now we assume that the components of $\omega$ are all distinct so that we can make a symplectic linear change of variables that diagonalize the quadratic part, hence
\[ H(z)=\omega.\mathcal{I}+V(z) \]
where $\mathcal{I}=\mathcal{I}(z)$ are the ``formal actions", that is
\[ \mathcal{I}(z)=(\mathcal{I}_1,\dots ,\mathcal{I}_n)=\frac{1}{2}(x_1^2+y_1^2, \dots, x_n^2+y_n^2)\in\R^n\ \text{
\rm and}\ \omega\!.\mathcal{I}=\displaystyle{\sum_{j=1}^{n}}\omega_j\mathcal{I}_j.
\]

 If all the components of $\omega$ have all the same sign, then $H$ is a
Lyapounov function and the fixed point is stable. But in a general setting,
one has to take into account the higher order terms $V(z)=O(|z|^3)$.

\medskip

 Given a solution $z(t)$, if we denote $\mathcal{I}(t)=\mathcal{I}(z(t))$, then $||\mathcal{I}(t)||$ is essentially the distance of $z(t)$ to the origin and $||\mathcal{I}(t)-\mathcal{I}(0)||$ measures the deviation of $z(t)$ from its original position.

 Without loss of generality, we can assume that we consider a Hamiltonian defined on the ball $B_R$ of radius $R<1$ in $\R^{2n}$. Then by analyticity, we extend the resulting Hamiltonian as a holomorphic function on some closed complex ball $D_s$ of radius $s<1$ in $\C^{2n}$. If we define $\mathcal{A}_s$ as the space of holomorphic Hamiltonians on $D_s$ which are real valued for real arguments, and $||.||_s$ its usual supremum norm, then we are led to consider
\begin{equation}\label{Ham1}
\begin{cases}
H(z)=\omega.\mathcal{I}+f(z) \\
H \in \mathcal{A}_s\ ;\ f(z)=O(z^3).
\end{cases}\tag{$A$}
\end{equation}

\medskip

To obtain results of stability, one of the main tool is to construct normal forms via classical averagings, and in this case these are the so-called {\it Birkhoff normal forms}.

 For an integer $m\geq 1$, assuming $\omega$ is non-resonant up to order $2m$, that~is
\[ k.\omega \neq 0, \quad k\in\Z^n, \; 0<||k||_1 =\sum_{j=1}^n | k_j |\leq 2m  \]
then there exists an analytic symplectic transformation $\Phi_m$ close to identity such that $H\circ\Phi_m$ is in Birkhoff normal form up to order $2m$, that is
\[ H\circ\Phi_m(z)=h_m(\mathcal{I})+f_m(z)\ \text{\rm where}\ h_m(\mathcal{I})=
\displaystyle{\sum_{k=1}^m}\mathcal{B}^{(k)}(\mathcal{I})\]
with homogenous polynomial $\mathcal{B}^{(k)}$ of degree $k$ in the $\mathcal{I}$ variables and the remainder $f_m$ of order $z^{2m}$ (see \cite{Bir66} or \cite{Dou88}, \cite{HZ94} for a more recent exposition). The polynomials $\left(\mathcal{B}^{(k)}\right)_{k\in\N^*}$ are uniquely defined (even though the transformations are not) and are usually called the {\it Birkhoff invariants}. The transformed Hamiltonian is therefore the sum of an integrable part $h_m$, for which the origin is trivially stable in the sense that $\mathcal{I}(t)$ is constant for all times, and a much smaller perturbation $f_m$ for a small enough neighbourhood of the origin. Moreover, if $\omega$ is non-resonant up to any order, we can even define a formal symplectic transformation $\Phi_\infty$ and a formal power series $h_\infty=\displaystyle{\sum_{k\geq 1}}\mathcal{B}^{(k)}$ such that

\begin{equation}
\label{hinfini}H\circ\Phi_\infty(z)=h_\infty(\mathcal{I}).
\end{equation}

However, for a certain topology on the coefficients, the transformation $\Phi_\infty$ is generically divergent as it was proved by Siegel \cite{Si41} and the convergence properties of  the series $h_\infty$ are even more subtle. Actually, Perez-Marco (\cite{PM03}) has
showed that if $h_\infty$ is divergent for some $H\in\mathcal{A}_s$, then the complete Birkhoff normal form $h_\infty$ is divergent for a "typical"
$H\in\mathcal{A}_s$. But nevertheless, we cannot ensure that $h_\infty$ is generically convergent or divergent.

\medskip

The first kind of stability result is given by an application of KAM theory. Assume that $\omega$ is non-resonant up to order $4$ so that our system reduce to
\[ H(z)=h_2(z)+f_2(z)=\omega.\mathcal{I}+\mathcal{B}^{(2)}(\mathcal{I})+f_2(z). \]
where $\mathcal{B}^{(2)}$ is a quadratic form.

 Under the so-called Kolmogorov non-degeneracy condition (or the twist condition in the context of an elliptic fixed point) : the quadratic form $\mathcal{B}^{(2)}$
is {\it non-degenerate}, KAM theory (see \cite{AKN97} and \cite{DG96} for the specific case of an elliptic fixed point) ensures stability for all times for most of the initial conditions. That is, for a large set of solution $z(t)$ :
if $|\mathcal{I}(0)|\leq\rho$ then the variations $|\mathcal{I}(t)-\mathcal{I}(0)|$ are of order $\rho$ for all times, hence we have Lyapounov stability. More precisely, KAM theory yields the existence of invariant Lagrangian tori which form a set of positive measure with a density going to one at the origin. However because of the dimensions, it is only for $n=2$ that we can deduce perpetual stability results. Actually, for $n=2$ the invariant tori are 2-dimensional and divide the 3-dimensional energy level, therefore the solutions of the perturbed system are global and bounded over infinite times. On the other
hand, an arbitrary large drift of the orbits can still occur for $n\geq
3$ in the complement of the invariant tori which is a dense set with a complicated
topology. In fact it is believed that, except for two degrees of freedom systems (and of course one degree), ``generic" elliptic fixed points are unstable (see \cite{DLC83} and \cite{Dou88} for examples and \cite{KMV04} for an announcement in a generic setting).

 Thus, for $n>2$, KAM theory provides results of stability only in the sense
of measure and theorems of stability which are valid for an open set of initial condition can only be proved over {\it finite} times.

\begin{definition}
The origin is stable over {\it finite times} if $|\mathcal{I}(0)|\leq\rho$ small enough implies
\[ |\mathcal{I}(t)-\mathcal{I}(0)|=R(\rho)\ \text{\rm for}\ |t|<T(\rho)\]
where $R(\rho )$ is of order $\rho$ and $T(\rho )$ is at least of order
$\rho^{-1}$ (otherwise we have a trivial bound).

\medskip

 We have {\bf polynomial stability} if $T(\rho )$ is at least of order
$\rho^{-m}$ for some fixed integer $m\in\N^*$.

 We have {\bf exponential stability} if $T(\rho )$ is at least of order
$\exp\left( C\rho^{-\alpha}\right)$ for some fixed constants $C>0$ and $\alpha
>0$.

 We have {\bf superexponential stability} if $\exp\left( C\rho^{-\alpha}\right)$
 is negligible with respect to $T(\rho )$ for {\it any} constants $C>0$ and $\alpha >0$.
\end{definition}

\medskip

 Using the mean value theorem, a polynomial time of stability can be easily obtained once a Birkhoff normal form up to a finite order has been built. The first result of this kind around an elliptic fixed point was obtained by Littlewood in 1959 (\cite{Lit59}).

 To obtain results of exponential stability, there is basically two distinct methods.

In the first approach, we call \text{\it Birkhoff type estimates}, one assumes that the frequency $\omega$ is a $(\gamma ,\tau )${\it -Diophantine vector} for some positive constants $\gamma$ and $\tau$, hence $\omega\in\Omega_{\gamma ,\tau}$ where :
\begin{equation}
\Omega_{\gamma ,\tau}=\left\{\omega\in\R^n\ {\rm such\ that}\ \vert
k.\omega\vert\geq{{\gamma}\over{\vert\vert k\vert \vert_\infty^\tau}}\ {\rm for\ all}\ k\in\Z^n\backslash\{ (0,\ldots ,0\}\right\}
\end{equation}

We recall that the measure of the complementary set of $\Omega_{\gamma ,\tau}$ is of order ${\cal O}(\gamma )$ for $\tau >n-1$.

 In particular $\omega$ is completely non-resonant, hence one can perform any finite number $m$ of Birkhoff normalizations, and since we have a control on the small divisors, one can estimate the size of the remainder $f_m$.

 Under these assumptions, one can prove (\cite{BG86}, \cite{GDFGS89}) that for a real-analytic system, the ``action variables'' $\mathcal{I}$ become quasi integrals over
exponentially long times, more specifically :

\begin{theoreme}

{\it Consider an elliptic fixed point for a real analytic Hamiltonian system of type $(A)$.

 There exists positive constants $C_1$, $C_2$, $C_3$, $C_4$ which depend
only of the analyticity width $s$, the size $\vert\vert H\vert\vert_s$ of
the Hamiltonian $H$, the number of degree of freedom $n$ and $\tau$ the exponent
in the Diophantine condition such that if $\rho <C_1\gamma$, an arbitrary solution $(x(t),y(t))$ with an initial action $\vert\vert\mathcal{I}(0)\vert\vert
\leq\rho$ is defined at least over an exponentially long time and satisfies~:
\begin{equation}
\left\vert\left\vert \mathcal{I}(t)-\mathcal{I}(0)\right\vert\right\vert\leq C_2\rho\ {\rm if}\ \vert t\vert\leq C_3\exp\left( C_4\left(\frac{\gamma}{\rho}
\right)^{{{1}\over{1+\tau}}}\right) .
\end{equation}
}
\end{theoreme}

\medskip

 The proof rely on the existence of a normalizing transformation up
to an exponentially large order $m$ comparable with the inverse of the distance
to the origin. This is possible since we have lower bounds on the small denominators and the growth of the coefficients in the normalizing expansion is reduced to a combinatorial problem (the "great multipliers" along Poincar\'e terminology).
Actually, under the assumptions of linear diophantine frequencies and the
real analyticity of considered Hamiltonian, there exists a decreasing sequence
of radii $s_m\sim \gamma /m^{1+\tau}$ for $m\in\N^*$ such that the normalizing transformation and the Birkhoff normal up to order $m$ can be build on the
ball of radius $s_m$ with a remainder of size comparable with $s_m^m.(m!)^{1+\tau}\!\sim\! (m!/m^m)^{1+\tau}$ which is dominated by $\exp\left( -(1+\tau )m\right)$ or $\exp\!\left(\!\! -C\left(\frac{\gamma}{s_m}\right)^{{{1}\over
{1+\tau}}}\!\!\right)$.

 Finally, since the averaged Hamiltonian is integrable, the speed of drift of the action variables is very slow.

 More specifically, there exists a sequence of annuli of decreasing radii $(s_m)_{m\in\N^*}$ such that :
$$s_{m+1}\!\leq\! \vert\vert\mathcal{I}(0)\vert\vert\! <\! s_m\!\!\Longrightarrow
\!\!\left\vert\left\vert\mathcal{I}(t)-\mathcal{I}(0)\right\vert\right\vert\!\leq\! C_2 s_m\ {\rm if}\ \vert t\vert\!\leq\! C_3\exp\left(\!\! C_4\left(\!\frac{\gamma}
{s_m}\!\right)^{{{1}\over{1+\tau}}}\!\!\right)$$

 In order to obtain uniform estimates with respect to the initial radius
 $\rho$, as in our theorem 1.2, it only remains to slightly lower the coefficients $C_2, C_3, C_4$. Treschev continuous averaging method (\cite{Tr02}) would allow to obtain directly a result with respect to the initial radius $\rho$ and to avoid this sequence $(s_m)_{m\in\N^*}$ but we will stay in the setting of Giorgilli et al (\cite{GDFGS89}, see also \cite{DG96}).

\vskip0,5truecm

The second approach, we call \text{\it Nekhoroshev type estimates}, is fundamentally different. One just assume that $\omega$ is non-resonant up to order $4$
so that the Hamiltonian reduces as for KAM theory to
\[ H(z)=h_2(z)+f_2(z)=\omega.\mathcal{I}+\mathcal{B}^{(2)}(\mathcal{I})+f_2(z). \]
where $\mathcal{B}^{(2)}$ is a quadratic form. But this time we require that the non-linear part $h_2(z)=\omega.\mathcal{I}+\mathcal{B}^{(2)}(\mathcal{I})$ is {\it convex} in actions $\mathcal{I}$, which is equivalent to $\mathcal{B}^{(2)}$ being sign definite. Under those assumptions, it was predicted by Lochak (\cite{Loc92}) and proved independently by Niederman (\cite{Nie98}) and Benettin, Fasso and Guzzo (\cite{BFG98}) that exponential stability holds. Their proofs are based on the implementation of Nekhoroshev estimates in cartesian coordinates, but they are radically different : the first one uses Lochak method of periodic averagings and simultaneous Diophantine approximations while the second is based on Nekhoroshev original mechanism. The proof of Niederman was later clarified by Poschel (\cite{Pos99a}).

Actually, Birkhoff estimates rely on the existence of an accurate integrable normal form while Nekhoroshev theory (see \cite{Nek77} or \cite{Nie11} for
a survey) gives exponential stability for all the solutions even if resonances of {\it low order} are present. This problem is overcome by a suitable transversality condition called {\it steepness} on the geometry of the integrable system. Steepness is a generic condition infinitely differentiable functions which is implied by (quasi)-convexity.

\bigskip

 A remarkable result of Morbidelli and Giorgilli (\cite{MG95}) clearly shows that the two previous results which come respectively from Hamiltonian perturbation
theory and Nekhoroshev's theorem are {\it independent} and can be {\it superimposed}.

\medskip

\begin{theoreme}

{\it Consider an elliptic fixed point for a real analytic Hamiltonian system of type $(A)$ with a Diophantine normal frequency $\omega$ and a sign definite
torsion, that is : $\mathcal{B}^{(2)}$ is a sign definite quadratic form.

 There exists positive constants $C_1$, $C_2$, $C_3$, $C_4$ which depend
only of the analyticity width $s$, the size $\vert\vert H\vert\vert_s$ of
the Hamiltonian $H$, the number of degree of freedom $n$ and $\tau$ the exponent
in the Diophantine condition such that if $\rho <C_1\gamma$, an arbitrary solution $(x(t),y(t))$ with an initial action $\vert\vert\mathcal{I}(0)\vert\vert
\leq\rho$ is defined at least over a superexponentially long time and satisfies~:
\begin{equation}
\left\vert\left\vert \mathcal{I}(t)-\mathcal{I}(0)\right\vert\right\vert\leq C_2\rho\ {\rm if}\ \vert t\vert\leq C_3\exp\left(\exp\left( C_4\left(\frac{\gamma}
{\rho}\right)^{{{1}\over{1+\tau}}}\right)\right)
\end{equation}
}
\end{theoreme}

 The proof starts with a Birkhoff's normal form with an exponentially small remainder like in the theorem 1.2. Then, Nekhoroshev's theory is applied with a perturbation which is already exponentially small, hence we obtain a superexponential time of stability. More specifically, we need the
(quasi-)convexity of the Birkhoff normal forms at all order. Since (quasi-)
convexity
is an open condition for the $\mathcal{C}^{(2)}-$topology, the convexity
of the second Birkhoff invariant $\mathcal{B}^{(2)}$ implies the same property
for $h_m$ ($\forall m\in\mathbb{N}^*$).

\medskip

 One of the main problem to extend the later kind of results in a generic setting is that we need a property on the Birkhoff invariants at all order
which is {\it not} an open condition.

 In this direction, Bounemoura (\cite{Bo11a}) has proved that these reasonings
remain valid in the case of a non resonant elliptic fixed point under a generic
condition on the associated Birkhoff normal form $h_\infty$, more specifically :

\begin{theoreme}

{\it Consider an elliptic fixed point for a real analytic Hamiltonian system of type $(A)$ with a non-resonant normal frequency $\omega$.

 If its Birkhoff normal form $h_\infty$ belongs to a generic set in the space of the formal power expansions $\mathbb{R}[[X]]$ (genericity is in the sense
of prevalence, see the next section) the estimates of superexponential stability given in the previous theorem remain valid.
}\end{theoreme}

%\medskip

\section{Setting and statement of our results :}

%\medskip

 Bounemoura gave only a partial answer to the problem of superexponential
stability of elliptic fixed position. Actually, we cannot ensure that generically the Birkhoff normal form associated to a non resonant elliptic fixed point belongs to the set introduced by Bounemoura even if it is a big set in the sense of measure. This problem is not straightforward since we have very few general information on the Birkhoff normal forms. Actually, in the context
of KAM theory, we are only aware of one result of Eliasson, Fayad and Krikorian
(so far unpublished) which involves {\it all} the Birkhoff invariants and ensure the existence of invariant tori which accumulate an elliptic fixed point in a Hamiltonian system with a Diophantine linear spectrum but without any torsion condition.

 Here, we provide a complete answer to our problem of stability.

\medskip

 We first specify our space of Hamiltonians
$$\mathcal{H}_s=\left\{ H(x,y)=\omega.\mathcal{I} +f(x,y)\ \text{\rm where}\
\left\vert\begin{split}
 H &\in \mathcal{A}_s,\\
 \omega &\in \R^n\ \text{\rm and}\ f(x,y)=O_3 (x,y) .
\end{split}\right.\right\}$$
we recall that $\mathcal{A}_s$ is the space of holomorphic Hamiltonians on $D_s\subset\mathbb{C}^{2n}$ the ball of radius $s$ around the origin which are real valued for real arguments equipped with the supremum norm $\vert\vert .\vert\vert_s$ on $D_s$, hence it is a {\it Banach} space.

\medskip

 Without loss of generality, we can assume that $s<1$.

\bigskip

 For an homogenous polynomial $P(\mathcal{I})=\!\!\!\displaystyle{\sum_{\vert l\vert =k, l\in\N^{2n}}}\!\! p_k \mathcal{I}^k$, as in \cite{GHK06} we consider the
following Euclidean norm (also called Bombieri norm)
$$\vert\vert P\vert\vert_k =\sum_{\vert l\vert =k}\frac{\vert p_l\vert^2}{C_k^l}$$
with the multinomial coefficient $C_k^l$.

 This norm is invariant by isometry in~$\C^{2n}$ (see \cite{GHK06}).

\bigskip

  There exists a Borelian measure canonically associated to the norm $\vert\vert P\vert\vert_k$ on the space $\mathcal{P}_k$ of homogenous polynomials of degree $k$.

 \medskip

 We consider the {Hilbert brick}
$$HB^\N =\left\{ h(\mathcal{I})=\sum_{k\in\N^*}h_k(\mathcal{I})\ \text{with}\ \vert\vert h_k\vert\vert_k\leq 1\ (\forall k)\right\}$$
it is a {\bf compact} set for the product topology when the space of polynomial
is identified as the infinite product $\displaystyle{\prod_{k\in\N^*}}\mathcal{P}_k$.

This compact can be continuously embedded in our Banach space of Hamiltonian
$\mathcal{H}_s$ with the assumption that $s<1$, hence it is a compact set in $(\mathcal{H}_s,\vert\vert .\vert\vert_s)$.

\medskip

 If we denote by $\mu_k$ the Borelian probability measure on the unit ball of $\mathcal{P}_k$
for all $k\in\N^*$, the embedding of $HB^\N$ in $(\mathcal{H}_s,\vert\vert .\vert\vert_s)$ can be equipped with the image of the product measure $\mu =\displaystyle{\prod_{k\in\N^*}}
\mu_k$ which define a {\it compactly supported Borelian probability measure} in $(\mathcal{H}_s,\vert\vert .\vert\vert_s)$.

 This setting was previously considered by Gorodetski, Hunt and Kaloshin
 (\cite{GHK06}).

\medskip

\begin{theoreme}{\bf (Main theorem)}
For an arbitrary $H\in\mathcal{H}_s$, the modified Hamiltonian $H+h$ is {super-exponentially stable} for almost all $h\in HB^\N$ with respect to the probability $\mu$.

More precisely, there exist positive constants $a$, $C_1$, $C_2$, $C_3$ such that for a small enough $\rho$, every solution of $H+h$ with initial condition $|\mathcal{I}(0)|<\rho$ satisfies
\[ |\mathcal{I}(t)-\mathcal{I}(0)|< C_1\rho\ f\! or\ |t| < \exp\left( C_2\rho^{-C_3\vert
\ln(\rho )\vert^a}\right) .  \]
\end{theoreme}

 In other words, the set of super exponentially stable Hamiltonians is prevalent in our Banach space since $\mu$ defines a compactly supported Borelian measure which is transverse to the set of non super-exponentially stable Hamiltonian (see \cite{OY05} and \cite{HK10}).

 We do not see obstructions to obtain similar results for Lagrangian tori
 (see the remarks at the end of this section).

\medskip

\section{Scheme of the proof :}

\medskip

 Either in this paper or in Bounemoura work (\cite{Bo11a}), the proofs rely
crucially on the extension of Nekhoroshev theory obtained in \cite{Nie07}
for a much wider class of unperturbed Hamiltonian than the steep integrable Hamiltonian considered originally in \cite{Nek77}.

 For instance, it can be proved (\cite{Bo11a}) that almost all quadratic
form in action variables fall into this class of Diophantine Morse
functions while only sign definite quadratic forms are steep according to
Nekhoroshev definition.

 Actually, the definition of Diophantine Morse functions involves countable quantitative transversality conditions which are stated in adapted coordinates. It is inspired on one hand by the steepness condition introduced by Nekhoroshev (\cite{Nek77}) and on the other hand by the quantitative Morse-Sard theory of Yomdin (\cite{Yom83} and \cite{YC}) where we consider "nearly-critical" points which are quantitatively non degenerate.

 The application of Yomdin quantitative Morse Sard theory allows to prove
(\cite{Nie07}) that for an arbitrary real analytic integrable Hamiltonian $h$ defined in a neighbourhood of the the origin in the action space $\R^n$, then for almost any parameter $\omega\in\R^n$, the modified integrable Hamiltonian $h_\omega (x)=h(I)-\omega .I$ is exponentially stable (i.e. : Nekhoroshev estimates are valid).

 In the case of an elliptic equilibrium, we have to deal with cartesian coordinates instead of action-angles ones and we cannot use the latter since they become singular at the origin, this yields specific problems detailed in \cite{Loc92}
and overcomes in \cite{BFG98}, \cite{Nie98} for an elliptic point with a
sign definite torsion. In \cite{BN12}, Nekhoroshev estimates are
proved in a generic setting by means of successive periodic averagings and
this method can be implemented in cartesian coordinates. These observations
allowed Bounemoura \cite{Bo11a} to prove that a Morse Diophantine condition is satisfied at all order on a full Lebesgue measure set of Birkhoff invariants
which is strong enough to prove superexponential results of stability over
times of order $\exp (\exp (C/r^a ))$.

 As it will be specified in the sequel, the set considered by Bounemoura
 is not big enough to ensure that generically the Birkhoff normal form of
 a non resonant elliptic fixed point falls into this good set. We will
 have to consider more refined estimates on the set of bad parameter $\omega\in\R^n$
 according to the following proposition which is a direct extension of (\cite{Nie07}, theorem 3.2.6) or (\cite{BN12}, theorem 2.2) :
 \begin{proposition} \label{estbadpara}
Consider an integrable Hamiltonian
$$\mathrm{H}(x,y)=h(\mathcal{I})=h\left(\frac{
x^2+y^2}{2}\right)$$
holomorphic on the ball of radius $s$ around the origin
in $\mathbb{C}^{2n}$.

 There exists an open set $\mathcal{C}^{(\rho )}_\varepsilon\subset
\mathcal{B}_\rho$ where $\mathcal{B}_\rho$ is the ball of radius $\rho$ around the origin in $\mathbb{R}^n$ such that
\begin{equation*}
{\rm Vol}(\mathcal{C}^{(\rho )}_\varepsilon )< C_\rho \varepsilon^\alpha
\end{equation*}
where $C_\rho >0$ depends only of $n$, $s$, $\rho$, $\vert\vert\mathrm{H}\vert\vert_s$ and $\alpha >0$ depends only of $n$.

 For $\omega\in\mathcal{B}_\rho\backslash\mathcal{C}^{(\rho )}_\varepsilon$, Nekhoroshev estimates of stability can be proved for the modified integrable Hamiltonian $\mathrm{H}_\omega (x,y)=h_\omega (\mathcal{I})=h(\mathcal{I})-\omega .\mathcal{I}$ with an holomorphic perturbation $f(x,y)$ of size $\vert\vert f\vert\vert_s <\varepsilon$.

\end{proposition}

\begin{remarque}\label{remsdm}{\rm
Note also that $h$ does not need to be polynomial, we only require that the
Hamiltonian is smooth enough.}
\end{remarque}

\medskip

 Now, following the strategy of Morbidelli and Giorgilli (\cite{MG95}), the estimates of superexponential stability given by Bounemoura rely
on the existence for all $m\in\mathbb{N}^*$ of a Birkhoff normalization $\Phi_m$
specified in theorem 1.2. which integrate $H$ up to order $2m$ on the ball of radius $s_m >0$ around the origin in $\mathbb{C}^n$, hence $H\circ\Phi_m(z)\!
=\! h_m(\mathcal{I})+f_m(z)$ where $h_m\in\mathbb{R}_m[X]$ and the remainder $f_m$ is of order $z^{2m}$. Moreover, the size of the perturbation $\vert\vert f_m\vert\vert_{s_m}$ is comparable with $\exp\left( -(1+\tau )m\right)\vert\vert f_1\vert\vert_{s_1}$ and $\vert\vert f_1\vert\vert_{s_1}\!\sim\! s_1$ (see \cite{GDFGS89}, \cite{DG96}).

 Then, in order to obtain superexponential estimates of stability, we should be able to apply Nekhoroshev theory (with the proof given in \cite{BN12}) on the Birkhoff normal forms $\left( h_m\right)_{m\in\mathbb{N}^*}$ at all order $m$ with a perturbation $f_m$ which decrease geometrically.

 Going back to the genericity problem, according to the previous proposition for an arbitrary $m\in\N^*$, Nekhoroshev estimates can be proved on the ball of radius $s_m >0$ for the perturbation $f_m$ of $h_m(\mathcal{I}) +\omega .\mathcal{I}$ if $\omega\in\mathcal{B}_\rho\backslash\mathcal{C}^{(\rho )}_m$ where $\mathcal{C}^{(\rho )}_m$ is an open set which contain the critical parameters in the ball $\mathcal{B}_\rho
\subset\mathbb{R}^n$ with a volume which decrease geometrically with $m$ hence ${\rm Vol}\left(\displaystyle{\bigcup_{m\in\mathbb{N}^*}}\mathcal{C}^{(\rho )}_m\right) <\displaystyle{\sum_{m\in \mathbb{N}^*}}{\rm Vol}\left(\mathcal{C}^{(\rho )}_m\right) \sim {\rm Vol}\left(\mathcal{C}^{(\rho )}_1\right) <\infty$.

 In our infinite-dimensional dimensional setting, we will say that a formal
 expansion $h_\infty\in\mathbb{R}[[X]]$ is {\it superexponentially stable}
 if there exists a sequence of radii $\left( s_m\right)_{m\in\N^*}$ such
 that Nekhoroshev estimates can be proved for all $m\in\N^*$ with the integrable
 part given by the truncation ${\rm Proj}_{\mathbb{R}_m[X]}(h_\infty )$ at
 order $m$ and a perturbation $f_m$ of size $\vert\vert f_m\vert\vert_{s_m}\sim
 s_1/e^{(1+\tau )m}$\!\! .

  We note that we have ${\rm Proj}_{\mathbb{R}_m[X]}(h_\infty +\ell_\omega )={\rm Proj}_{\mathbb{R}_m [X]}(h_\infty )+\ell_\omega$ where $\ell_\omega (X)=\omega .X$ for all $\omega\in\mathbb{R}^n$ hence, for an arbitrary formal expansion $h_\infty$ in $\mathbb{R}[[X]]$, the modified expansion $h_\infty (\mathcal{I})+\omega .\mathcal{I}$ is superexponentially stable if  $\omega\in\mathcal{B}_\rho\backslash{\displaystyle{\bigcup_{m\in\mathbb{N}^*}}}
 \mathcal{C}^{(\rho )}_m$.

 From these estimates Bounemoura (\cite{Bo11a}) proves that, for an arbitrary
formal expansion $h_\infty\in\mathbb{R}[[X]]$, superexponential estimates
of stability are valid for $h_\infty (\mathcal{I})-\omega .\mathcal{I}$
on a ball of small enough radius $\varepsilon$ around the origin if $\omega\in\mathcal{B}_\rho
\backslash\mathcal{C}^{(\rho )}_\varepsilon$ where $\mathcal{C}^{(\rho )}_\varepsilon
\!\! =\!\!\!\displaystyle{\bigcup_{m\in \mathbb{N}^*}}\mathcal{C}^{(\rho )}_m$ with ${\rm Vol}\left(\mathcal{C}^{(\rho )}_\varepsilon\right)$ of the order of $\varepsilon^\alpha$ where $\alpha$ given in the previous proposition depends only of $n$.

 This implies the generic property stated in Theorem 1.4.

 Especially, Bounemoura (\cite{Bo11a}) shows the existence of a foliation
of the space of formal series $\mathbb{R}[[X]]$ in $n$-dimensional probe
spaces directed by the space of linear forms where the properties required
to prove superexponential estimates are satisfied at almost every point on the leaves for the $n$-dimensional Lebesgue measure. The issue is that when we want to come back in the initial (non-normalized) coordinates, a change of the linear frequencies implies a change of all the Birkhoff invariants and all the normalizing transformations. Hence it is impossible to recover a foliation of the initial space of Hamiltonian $\mathcal{H}_s$ by $n$-parameter
families where superexponential stability can be proved for almost all parameters
with the $n$-dimensional Lebesgue measure. We must use more global estimates
with an {\it infinite} number of parameters and the important observations in this direction are the following.

\medskip

 To state our results we have considered the embedding of $HB^\mathbb{N}$
in the Banach space $\mathcal{H}_s$. Here, we first consider the Hilbert
brick $HB^\mathbb{N}$ embedded in the space of formal expansion $\mathbb{R}[[X]]$ identified with the product space ${\prod_{\N^*}}\mathcal{P}_k$
equipped with the product topology (we recall that $\mathcal{P}_k$ is the space of homogenous
polynomial of degree $k\in\N^*$).

 Actually, despite we are in an infinite dimensional space,
the central estimates concern polynomial integrable Hamiltonian (the truncated
Birkhoff normal forms) and now we revisit the study of Bounemoura for polynomial
Hamiltonian.  For all $m\in\mathbb{N}^*$, we identify the space $\mathbb{R}_m[X]$
of polynomial of degree bounded by $m$ with the product $\displaystyle{\prod_{1\leq k\leq m}}\mathcal{P}_k$ of the spaces $(\mathcal{P}_k)_{k\in\mathbb{N}^*}$ of homogenous polynomial of degree $k\leq m$ equipped with the product topology.
Moreover, for all $m\in\N^*$, we denote the projection ${\rm Proj}_{\mathbb{R}_m[X]}
\left( HB^\mathbb{N}\right)$ by $HB^{\leq m}\subset\mathbb{R}_m[X]$ equipped
with the product measure $\mu_{\leq m}=\!\!\displaystyle{\prod_{1\leq k\leq m}}\!\!\mu_k$.

 We note that, for all $m\in\N^*$, we have $\mu =\mu_{\leq
m}\times\displaystyle{\!\!\prod_{k\geq m+1}}\!\!\mu_k$ and especially $\mu =\mu_1\times
\displaystyle{\prod_{k\geq 2}}\mu_k$ where $\mu_1$ is the classical Euclidean volume up to a scaling (to define a probability) on the unit ball in the space of linear forms.

\medskip

 As in the previous reasonings, for all $m\in\mathbb{N}^*$, we consider
an open set $X_\varepsilon^{(m)}$ in $HB^{\leq m}\subset\mathbb{R}_m[X]$ given by Yomdin Morse-Sard theory which include the expansions in $HB^{\leq
m}$ (considered as integrable Hamiltonians) where Nekhoroshev theory cannot be applied to prove superexponential results of stability on a ball of radius
$\varepsilon$. Using Bounemoura construction, we can prove that $X_\varepsilon^{(m)}$ intersects any affine subspace in $\mathbb{R}_m[X]$ directed by the vectorial space of linear forms along a set of volume dominated by $\left(\varepsilon /\xi^m\right)^\alpha$ with $\alpha >0$, given in proposition 3.1, which depends only of the number of degree of freedom $n$ and a constant $\xi >1$.

 With the fact that $\mu_{\leq m}=\mu_1\times\displaystyle{\prod_{k=2}^m}\mu_k$
where $\mu_1$ is the Euclidean volume and $\displaystyle{\prod_{k=2}^m}\mu_k$
is a probability measure, the application of Fubini theorem implies that the measure $\mu_{\leq m}\left( X_\varepsilon^{(m)}\right)$ remains is of the order of $\left(\displaystyle{\frac{\varepsilon}{\xi^m}}\right)^\alpha$.

 Now we come back to our initial functional space $\mathcal{H}_s$, for an arbitrary Hamiltonian $H\in\mathcal{H}_s$ we consider the application :
$$\begin{array}{ccll}
\mathcal{B}\mathcal{N}\!\mathcal{F}_m : &\mathcal{M}_m &\longrightarrow &\ \ \ \ \mathbb{R}_m [X]\ \ \ \ \ \ \ {\mathrm w}{\mathrm i}{\mathrm t}{\mathrm h}\ \ \ \mathbb{R}_m[X]\simeq\!\!\!\!
\displaystyle{\prod_{1\leq k\leq m}}\mathcal{P}_k\\
& \left( P_1, \ldots , P_m\right) &\longmapsto &\left( Q_1, \ldots , Q_m\right)
\end{array}$$
where $Q_k$ is the Birkhoff invariant of order $k\in\{ 1,\ldots ,m\}$ of $H+P$ for $P={\displaystyle\sum_{k=1}^m} P_k$ and $\mathcal{M}_m$ is an open set in $\mathbb{R}_m [X]$ where the Birkhoff normal form at order $m$ can be build.

 We have :

\medskip

\begin{proposition}

 For $k=1$, the transformation is $Q_1=P_1+\mathcal{B}^{(1)}(H)$ where $\mathcal{B}^{(1)}(H)$ is the first Birkhoff invariant of $H$ hence it is a translation.

 For $k\in\{ 2,\ldots ,m\}$ :
$$Q_k=P_k +\mathcal{F}_k (P_1,\ldots ,P_{k-1})$$
where $\mathcal{F}_k$ is a differentiable function over $\mathcal{M}_m$, hence the mapping $\mathcal{B}\mathcal{N}\!\mathcal{F}_m$ is one-to-one with a Jacobian equal to 1.

\end{proposition}

 A similar property on the Newton interpolation polynomials was used by
Gorodetski-Hunt-Kaloshin (\cite{GHK06}).

\medskip

 Despite our transformation has a unit Jacobian, it does not preserve the
measure $\mu_{\leq m}$ since it does not leave invariant the projection $HB^{\leq m}$ of the Hilbert Brick~$HB^\mathbb{N}$.

 Actually, for an analytic Hamiltonian, the Birkhoff normal form should admit a Gevrey regularity and not remain analytic (we recall that we have no results on the generic divergence or convergence of the Birkhoff normal forms). Hence, the Birkhoff invariants should growth much more rapidly then the coefficient of the initial Hamiltonian and the set $HB^{\leq m}\cap\mathcal{M}_m$ is sent to a much bigger domain by the mapping $\mathcal{B}\mathcal{N}\!\mathcal{F}_m$.
Consequently, we cannot compute the measure $\mu_{\leq m}$ of an image set in $\mathcal{B}\mathcal{N}\!\mathcal{F}_m\left( HB^{\leq m}\cap\mathcal{M}_m\right)$ by the change of variables formula.

\medskip

 In order to overcome this problem, we rescale the normalized variables in $\mathbb{C}^{2n}$ on the ball of radius $s_m$ around the origin where the Birkhoff transformation $\Phi_m$ at order $m$ is defined. Denoting the initial
variables by $z\in\mathbb{C}^{2n}$, we consider $z'\in\mathbb{C}^{2n}$ defined by $z=\Phi_m (s_m z')$ which correspond to a conformal symplectic transformation giving the Hamiltonian~:
\begin{equation}
K_m(z')=\frac{1}{s_m^2}H\circ\Phi_m (s_m z')\ \text{\rm or}\ K_m(x',y')=\frac{1}
{s_m^2}H\circ\Phi_m (s_m x',s_m y').
\end{equation}

 In these rescaled variables, the integrable part $k_m(z')=\frac{1}{s_m^2}h_m(s_m^2 \mathcal{I}(x',y'))$ belongs to $HB^{\leq m}$.

 Our initial issue is settled but the price to pay is a huge deformation of the measure $\mu_{\leq m}$ since the Jacobian of the transformation $\widetilde
{\mathcal{B}\mathcal{N}\!\mathcal{F}}_m$ which maps to the rescaled Birkhoff invariants (denoted $\left({\tilde Q}_1,\ldots ,{\tilde Q}_m\right)$ in the
sequel) is now $s_m^{-{\rm D}(m)}$ where ${\rm D}(m)$ growth like dimension of the space of polynomial $\mathbb{R}_m[X]$ of degree $m$ in $n$ variables hence ${\rm D}(m)$ is of the order of~$m^n$.

\medskip

 Let ${\widetilde X}_\varepsilon^{(m)}$ be a measurable set given by Yomdin quantitative Morse Sard theory of polynomial integrable Hamiltonians in $HB^{\leq m}$ which include the set of critical Hamiltonians at order $m$ where Nekhoroshev
theory cannot be applied in the rescaled variables to prove superexponential
estimates of stability.

 We must bound the measure of the set of Hamiltonian in $HB^{\leq m}\cap\mathcal{M}_m$ which are sent in ${\widetilde X}_\varepsilon^{(m)}$ by $\widetilde{\mathcal{B}\mathcal
{N}\!\mathcal{F}}_m$, hence :
$${\displaystyle\int\limits_{\mathcal{M}_m}}\!\!\!\! 1_{{\widetilde X}_\varepsilon^{(m)}}\!\!
\left(\widetilde{\mathcal{B}\mathcal{N}\!\mathcal{F}}_m (P_1,\ldots ,P_m)\right) \!\! d\mu_{\leq m}(P)\! =\!\frac{1}{s_m^{{\rm D}(m)}}\!\displaystyle{\int\limits_{\widetilde
{\mathcal{B}\mathcal{N}\!\mathcal{F}}_m (\mathcal{M}_m)}}\!\!\!\!\!\!\!\!\!\!\!\! \!\! 1_{{\widetilde X}_\varepsilon^{(m)}}\!\!\left({\tilde Q}_1,\ldots ,{\tilde Q}_m\right)\!\! d\mu_{\leq m} ({\tilde Q})$$
and this is bounded by ${\displaystyle\frac{\mu_{\leq m}\left({\widetilde X}_\varepsilon^{(m)}\right)}{s_m^{{\rm D}(m)}}}$.

\medskip

 With the growth of ${\rm D}(m)$, the divisor $s_m^{-{\rm D}(m)}$ decreases
 in a {\it superexponential} way with respect to the order of normalization
 $m$. At this step, the set introduced by Bounemoura is not big enough to
 ensure that its preimage in the construction of the Birkhoff normal forms
 is a generic set in $\mathcal{H}_s$. More specifically, we cannot afford a {\it geometric} decrease of the volume the critical parameters as in the convex case or in Bounemoura construction with the set $X_\varepsilon^{(m)}$.

 Here, we need a superexponential decrease of the volume of ${\widetilde X}_\varepsilon^{(m)}$\!\!\! .

 In order to settle this issue, we will not use Nekhoroshev theory on the
 whole domain where the Birkhoff normalization at order $m$ is defined (hence
 on the ball of radius $s_m$) but on a much smaller domain which is the ball
 ${\widetilde{\mathcal{B}}}_{r_m}$ of radius $r_m$ around the origin in the
 rescaled variables (hence $r_m s_m$ in the normalized variables).

 After a normalization at order $m\in\N^*$ and our rescaling, the size of the perturbation is of the order of $(m!)^{1+\tau}(s_m r_m)^m\sim \exp
(-(1+\tau )m)r_m^m$ (the prefactor comes from the growth of the coefficients).
Hence the size of the perturbation is dominated by $r_m^m$ since the expansion
of the perturbation starts at order $m$ (this is another important point),
with the proposition 3.1 and the previous reasonings we have $\mu_{\leq m}
\left({\widetilde X}_\varepsilon^{(m)}\right)$ of the order of $(r_m)^{m\alpha}$\!\!\!
.

\medskip

 Now we come back in our infinite-dimensional setting and, for all $m\in\N^*$,
 we consider the set ${\mathbf{{\widetilde X}_{\mathbf\varepsilon}^{\mathbf(m)}}}\!\!\!
\subset\mathbb{R}[[X]]$ such that $h_\infty\in{\mathbf{{\widetilde X}_{\mathbf
\varepsilon}^{\mathbf(m)}}}$ if and only if ${\rm Proj}_{\mathbb{R}_m[X]}\left( h_\infty\right)$ is in $\left(\widetilde{\mathcal{B}\mathcal{N}\!\mathcal{F}}_m\right)^{-1}\!\! \left({\widetilde X}_\varepsilon^{(m)}\right)$ or $\widetilde{\mathcal{B}\mathcal{N}\!
\mathcal{F}}_m\!\!\left({\rm Proj}_{\mathbb{R}_m[X]}\left( h_\infty\right)\right)\!\!\in
\!\!{\widetilde X}_\varepsilon^{(m)}$.

\medskip

 For all $m\in\N^*$, the product $\displaystyle{\prod_{k\geq m+1}}\mu_k$ is a probability measure and we have $\mu =\mu_{\leq m}\times\!\!\!\displaystyle
{\prod_{k\geq m+1}}\mu_k$, hence : $\mu\left({\mathbf{{\widetilde X}_{\mathbf\varepsilon}^{\mathbf(m)}}}
\right)\!\! =\!\mu_{\leq m}\left(\left(\widetilde{\mathcal{B}
\mathcal{N}\!\mathcal{F}}_m\right)^{-1}\!\! \left({\widetilde X}_\varepsilon^{(m)}
\right)\right)$.

 With the previous reasonings for all $m\in\N^*$, we have the upper bounds~:
$$\mu\left({\mathbf{{\widetilde X}_{\mathbf\varepsilon}^{\mathbf(m)}}}\right)\!\!
\leq\frac{\mu_{\leq m}\left({\widetilde X}_\varepsilon^{(m)}\right)}{s_m^{{\rm D}(m)}}\sim\frac{(r_m)^{m\alpha}}{s_m^{{\rm D}(m)}}.$$

 Finally, the set of critical Hamiltonian $\mathbf{C}$ where Nekhoroshev
theory cannot be applied at some order $m\in\mathbb{N}^*$ to prove superexponential
results of stability is included in $\displaystyle{\bigcup_{m\in\mathbb{N}^*}}{\mathbf
{{\widetilde X}_{\mathbf\varepsilon}^{\mathbf(m)}}}$ which is of small measure if the sum
$${\displaystyle\sum_{m\in\mathbb{N}^*}\frac
{(r_m)^{m\alpha}}{s_m^{{\rm D}(m)}}}<\infty.$$

 The second radius $r_m$ is a {\it free} parameter which has to be chosen small enough to ensure the convergence of ${\displaystyle\sum_{m\in\mathbb{N}^*}\frac
{(r_m)^{m\alpha}}{s_m^{{\rm D}(m)}}}$.

 This imposes to choose a radius $r_m$ which is {\it exponentially small} with respect to the order of normalization $m$, here $m$ is of the order of $\log (r_m)^a$ where $a>0$ depends only of $n$. We cannot recover an exponentially small remainder with $m$ of the order of $1/r_m$ as in the convex case but only a remainder which is smaller than any polynomially small remainder with respect to the radius. This is an important point which comes exclusively from the expected divergence of the Birkhoff normal form (actually, as previously said, we know that this normal form is always obtained as a Gevrey divergent expansion in the Diophantine case but we don't know if it diverges or converges generically). Another misleading fact is that on the ball of radius $s_m r_m$, we could obtain directly a Birkhoff normal form which integrate the system up to a superexponentially small remainder with respect to $m$. But in this case, the Birkhoff normal form would give only a remainder which is exponentially small with respect to the radius while here we start with a perturbation which is smaller than any polynomially small remainder with respect to the radius before going to the exponential with Nekhoroshev estimates.

 Summarizing we have proved that there exists positive constants $a$, $C_1$,
 $C_2$, $C_3$ and a sequence of annuli of decreasing radii $(r_m)_{m\in\N^*}$ such that, for an arbitrary function $H\in\mathcal{H}_s$ and for almost all $h\in HB^\N$ with respect to the probability $\mu$, every solution of the Hamiltonian system governed by $H+h$ with an initial condition $r_{m+1}\leq
\vert\vert{\widetilde{\mathcal{I}}}(0)\vert\vert <r_m$ in the rescaled variables satisfies :
\[ \left\vert\left\vert{\widetilde{\mathcal{I}}}(t)-{\widetilde{\mathcal{I}}}(0)\right\vert\right\vert< C_1 r_m\ f\! or\ |t| < \exp\left( C_2 r_m^{-C_3\vert\ln( r_m )\vert^a}\right)\!\! .  \]

 In order to obtain uniform estimates with respect to the initial radius
 $\rho$, as in our theorem 1.4, it only remains to slightly lower the coefficients $C_1, C_2, C_3$ and the exponent $a$.

\medskip

\section{Comments and prospects :}

\medskip

 A natural extension of this work would be for the case of a Lagrangian invariant
torus for an Hamiltonian system. In this setting, we consider a system governed
by the Hamiltonian $\omega .I +f(I,\varphi )I^2$ where $f$ is analytic in
action-angle variables $(I,\varphi )\in\mathbb{R}^n\times\mathbb{T}^n$. The only missing ingredient to make the previous reasonings is a quantitative
construction of the Birkhoff normal form to a given order in action. Up to the author knowledge, the available constructions (\cite{Fas90}, \cite{PW94},
\cite{Pos93}) in the literature yield an integrable Birkhoff normal form with a perturbation which is small but with terms of all order in action
including low order terms. Here, we need a Birkhoff normal form of the kind
$h_m(I)+f_m(I,\varphi )I^m$ for $m\in\mathbb{N}^*$ to follow the scheme of
our proof. This is a purely technical issue and there is certainly no obstruction to overcome this problem.

  The second extension would be for the case of an Hamiltonian with a lower regularity, that is for Gevrey or differentiable Hamiltonians. With the classical
papers on Birkhoff normal forms in the differentiable case or the work of Popov-Mitev \cite{MP10} in the Gevrey case together with the article of Bounemoura \cite{Bo11b} on generic Nekhoroshev estimates for Gevrey or differentiable Hamiltonian we have all the required technical tools.

 A final interesting extension would be for the case of invariant tori provided
by KAM theory. A model statement would be : for an analytic integrable hamiltonian which satisfy Kolmogorov (or Russman) non-degeneracy condition (\cite{AKN97}), KAM invariant tori are superexponentially stable for a small enough generic perturbation. The papers of Popov (\cite{Pop00}, \cite{Pop04}) on the existence of a global Birkhoff normal form for all the tori given by KAM theory should be useful in this setting.

\bigskip

\addcontentsline{toc}{section}{References}
\bibliographystyle{amsalpha}
%\bibliography{SuperExpEffArticle}

\newcommand{\etalchar}[1]{$^{#1}$}
\providecommand{\bysame}{\leavevmode\hbox to3em{\hrulefill}\thinspace}
\providecommand{\MR}{\relax\ifhmode\unskip\space\fi MR }
\providecommand{\MRhref}[2]{%
  \href{http://www.ams.org/mathscinet-getitem?mr=#1}{#2}
}
\providecommand{\href}[2]{#2}

\end{document}